\documentclass[11pt]{article}
\begin{document}
\renewcommand{\thefootnote}{\fnsymbol{footnote}}
\pagestyle{empty}
\setcounter{page}{1}
\def\beq{\begin{equation}}
\def\eeq{\end{equation}}
\def\bea{\begin{eqnarray}}
\def\eea{\end{eqnarray}}
\def\nn{\nonumber}
\def\cUh{{\cal U}_{\sf h}(osp(2|1))}
\def\Uh{{\bf U_h}(osp(2|1))}
\def\Uq{{\cal U}_{q}(osp(2|1))}
\def\U{{\cal U}(osp(2|1))}
\renewcommand{\thesection}{\arabic{section}}
\renewcommand{\theequation}{\thesection.\arabic{equation}}
\newfont{\twelvemsb}{msbm10 scaled\magstep1}
\newfont{\eightmsb}{msbm8} \newfont{\sixmsb}{msbm6} \newfam\msbfam
\textfont\msbfam=\twelvemsb \scriptfont\msbfam=\eightmsb
\scriptscriptfont\msbfam=\sixmsb \catcode`\@=11
\def\Bbb{\ifmmode\let\next\Bbb@\else \def\next{\errmessage{Use
      \string\Bbb\space only in math mode}}\fi\next}
\def\Bbb@#1{{\Bbb@@{#1}}} \def\Bbb@@#1{\fam\msbfam#1}
\newfont{\twelvegoth}{eufm10 scaled\magstep1}
\newfont{\tengoth}{eufm10} \newfont{\eightgoth}{eufm8}
\newfont{\sixgoth}{eufm6} \newfam\gothfam
\textfont\gothfam=\twelvegoth \scriptfont\gothfam=\eightgoth
\scriptscriptfont\gothfam=\sixgoth \def\frak{\frak@}
\def\frak@#1{{\fam\gothfam{{#1}}}} \def\frak@@#1{\fam\gothfam#1}
\catcode`@=12
%
%
%
\def\CC{{\Bbb C}}
\def\NN{{\Bbb N}}
\def\QQ{{\Bbb Q}}
\def\RR{{\Bbb R}}
\def\ZZ{{\Bbb Z}}
\def\cA{{\cal A}}          \def\cB{{\cal B}}          \def\cC{{\cal C}}
\def\cD{{\cal D}}          \def\cE{{\cal E}}          \def\cF{{\cal F}}
\def\cG{{\cal G}}          \def\cH{{\cal H}}          \def\cI{{\cal I}}
\def\cJ{{\cal J}}          \def\cK{{\cal K}}          \def\cL{{\cal L}} 
\def\cM{{\cal M}}          \def\cN{{\cal N}}          \def\cO{{\cal O}}
\def\cP{{\cal P}}          \def\cQ{{\cal Q}}          \def\cR{{\cal R}} 
\def\cS{{\cal S}}          \def\cT{{\cal T}}          \def\cU{{\cal U}}
\def\cV{{\cal V}}          \def\cW{{\cal W}}          \def\cX{{\cal X}}
\def\cY{{\cal Y}}          \def\cZ{{\cal Z}}
\def\qed{\hfill \rule{5pt}{5pt}}
\def\arcsinh{\mathop{\rm arcsinh}\nolimits}
\newtheorem{theorem}{Theorem}
\newtheorem{prop}{Proposition}
\newtheorem{conj}{Conjecture}
\newenvironment{result}{\vspace{.2cm} \em}{\vspace{.2cm}}

\pagestyle{plain}

\begin{center}

$\phantom{xx}$

\smallskip
\smallskip

{\LARGE {\bf ON NONSTANDARD QUANTIZATIONS OF $osp(2|1)$ SUPERALGEBRA 
VIA CONTRACTION AND MAPPING}}
\\[0.5cm]

{B. ABDESSELAM$^{a,b}$\footnote{E-mail address: boucif@yahoo.fr}, 
R. CHAKRABARTI$^{c,}$\footnote{E-mail address: ranabir@imsc.res.in}, 
A. HAZZAB$^{d}$ 
and 
A. YANALLAH$^{b}$}

\smallskip
\smallskip 

{\it $^{a}$Laboratoire de Physique Th\'eorique d'Oran (LPTO), 
Université d'Oran Es-Sénia, 31100-Oran, Alg\'erie\\
$^{b,}$\footnote{Laboratoire de recherche agrée par le MESRS dans le 
cadre du fond national de la recherche et du developpement 
technologique.}Laboratoire de Physique Quantique de la 
Matière et Modélisations Mathématiques (LPQ3M),\\ 
Centre Universitaire de Mascara, 29000-Mascara, Alg\'erie \\ 
$^{c}$Department of Theoretical Physics, University of Madras, Guindy 
Campus, Madras 600025, India\\
$^{d}$Laboratoire de Modlisation et Mthode de Calcul, Centre 
Universitaire de Saida, 20000-Saida, Alg\'erie}. 
\end{center}

$\underline
{\phantom{XXXXXXXXXXXXXXXXXXXXXXXXXXXXXXXXXXXXXXXXXXXXXXXXX}}$

\begin{abstract}

\noindent We develop a generic representation-independent 
contraction procedure for obtaining, for instance, $R_{\sf h}$ 
and $L$ operators of arbitrary dimensions for the quantized
$\cUh$ algebra corresponding to the classical $r_2$ matrix
from the pertinent quantities of the standard $q$-deformed
$\Uq$ algebra. Also the quantized $\Uh$ algebra corresponding 
to the classical $r_1$ matrix comprising of the generators of the
classical $sl(2)$ algebra is obtained in terms of a nonlinear basis 
set.      

\noindent $\underline
{\phantom{XXXXXXXXXXXXXXXXXXXXXXXXXXXXXXXXXXXXXXXXXXXXXXXXX}}$

\end{abstract}

\setcounter{equation}{0}
\section{Introduction}
Quantum deformations of the Lie superalgebra $osp(2|1)$ have been 
studied  \cite{KR89}-\cite{BLT03} extensively both from the point 
of view of investigating integrable physical models, and also 
because of their intrinsic mathematical importance. It has been 
recently demonstrated \cite{JS98} that three distinct bialgebra 
structures exist on the classical $osp(2|1)$ superalgebra, and all of 
them are coboundary. The classical Lie superalgebra $osp(2|1)$ has 
three even $(h, b_{\pm})$ and two odd $(e, f)$ generators, obeying 
the following commutation relations: 
\bea
&&[h,\;e] = e,\qquad [h,\;f] = -f,\qquad \{e,\;f\}=-h,
\nn\\ 
&&[h,\;b_{\pm}] = {\pm}2 b_{\pm},\qquad 
\;[b_+,\;b_-] = h,\nn\\
&&[b_+,\;f]=e, 
\qquad [b_-,\;e]=f,\qquad b_+ = e^2,\qquad b_- = - f^2.
\label{eq:clalg}
\eea
The generators $(h, b_{\pm})$ form a subalgebra 
$sl(2) \subset osp(2|1)$. The classical universal enveloping algebra 
$\U$ is generated by the elements $(h, e, f)$. The irreducible 
representations are parametrized by an integer or a half-integer $j$, 
and are of dimension $(4j + 1)$.  For later use we quote here
a representation of the algebra (\ref{eq:clalg}):
\bea
h\,|j\;m> &=& 2m\,|j\;m>,\qquad 
e\,|j\;m> = |j\;m+1/2>,\nn\\
f\,|j\;m> &=& - (j + m)\,|j\;m-1/2>\quad \hbox{for}\;\;j - m\;\;
\hbox {integer},\nn\\   
\phantom{f\,|j\;m> }&=& (j - m + 1/2)\,|j\;m-1/2>\quad \hbox{for}
\;\;j - m\;\;\hbox {half-integer},
\label{eq:clrep}
\eea
where $m = j, j - 1/2,\dots,-(j - 1/2), -j$. The inequivalent 
classical $r$-matrices listed in Ref.\cite{JS98} read as follows: 
\bea
&&r_1 = h \wedge b_+,\qquad
r_2 = h \wedge b_+ - e \wedge e,\nn\\
&&r_3 = t (h \wedge b_+ + h \wedge b_- - e \wedge e - f \wedge f).
\label{eq:clrma}
\eea
The matrices $r_1$ and $r_2$ satisfy the classical Yang-Baxter 
equation, whereas the $r_3$ matrix satisfy the modified classical 
Yang-Baxter equation. The standard $q$-deformation of the $osp(2|1)$ 
algebra considered in Refs. \cite{KR89} and \cite{S90} corresponds 
to $r_3$, and it is a quasi-triangular algebra. The parameter 
$t$ in $r_3$ becomes irrelevant in quantization as it can be 
absorbed into the deformation parameter. Various applications 
of the algebra have been studied earlier. The $r_1$ matrix is comprised 
of the elements of the $sl(2)$ subalgebra. This allows quantization of 
the $osp(2|1)$ algebra using the inclusion $sl(2) \subset osp(2|1)$.
This is done \cite{CK98} by applying Drinfeld twist for the $sl(2)$
subalgebra to the full $osp(2|1)$ algebra. The Hopf superalgebra 
$\Uh$ obtained thereby is triangular. The quantization of the $r_1$ 
matrix has been obtained in Ref. \cite{CK98} in terms of the classical 
basis set with undeformed commutation relations but with coproduct 
structures deformed in a complicated manner.  
 
\par

The classical matrix $r_2$ has recently been quantized \cite{ACS03} 
using nonlinear basis elements. A class of invertible maps relating 
the quantized $\cUh$ algebra and its classical analog $\U$ has been 
found. The twist elements corresponding to these maps and the 
resultant universal ${\cal R}$ matrix have been evaluated as 
series expansions in the deformation parameter. In another work 
\cite{BLT03} the twist operator has been obtained in a closed-form 
expression in terms of the elements of the undeformed classical 
$osp(2|1)$ algebra. The universal ${\cal R}$ matrix obtained 
from this twist operator satisfies \cite{BLT03} the quantum 
triangularity condition. Another important issue observed \cite{K98} 
in this context is that the relevant quantum $R_{\sf h}$ matrix in the 
fundamental representation may be obtained {\it via} a contraction 
mechanism from the corresponding $R_q$ matrix of the standard $\Uq$ 
algebra in the $q \rightarrow 1$ limit. A generalization of this 
contraction procedure for arbitrary representations, though clearly 
desirable as it will allow us to systematically obtain various 
quantities of interest of the $\cUh$ algebra from the corresponding 
quantities of the $q$-deformed $\Uq$ algebra, has not been achieved 
so far.

\par

In the present work we complete the two above mentioned tasks. A generic
technic developed earlier \cite{ACC98}-\cite{ACC02} allowed us to 
extract quantum universal ${\cal R}$ matrix and ${\cal T}$ matrices for 
arbitrary representations for Jordanian ${\cal U}_{\sf h}(sl(N))$ 
algebra from the corresponding operators of the standard $q$-deformed 
${\cal U}_q(sl(N))$ algebra. A suitable adaptation of this methodology 
is used here to obtain the $L$ operator of the Jordanian $\cUh$ algebra 
obtained by quantizing the $r_2$ matrix. The standard FRT \cite{FRT90} 
procedure is then utilized, in conjunction with the $L$ operator 
obtained here, to derive the Hopf structure of the Borel subalgebra of 
the Jordanian $\cUh$ algebra. Our method may be readily used to derive 
$L$ operators for the higher representations and the corresponding 
${\cal T}$ matrices of the $\cUh$ algebra.

\par

In another problem considered here, we express the Jordanian deformed 
Hopf algebra $\Uh$ corresponding to the classical $r_1$ matrix in a 
{\it nonlinear} basis. Here we follow the approach in Ref. \cite{O92}, 
where the Jordanian ${\cal U}_{\sf h}(sl(2))$ algebra has been 
introduced in terms of a nonlinear basis set, while retaining the 
coproduct structure of these basis elements simple. 
A consequence of our choice of nonlinear basis elements is that, 
Ohn's ${\cal U}_{\sf h}(sl(2))$ algebra \cite{O92} explicitly
arises as a Hopf subalgebra in our $\Uh$ algebra. This feature is not
directly evident in the construction given in Ref. \cite{CK98}.
Moreover, while the algebraic relations of our $\Uh$ algebra 
are deformed, {\it the coproduct structures are considerably simple}. 
Our approach may be of consequence in building physical 
models of many-body systems employing coalgebra symmetry \cite{BH99}. 
Invertible nonlinear maps, and twist operators pertaining to these 
maps, exist connecting the deformed and the undeformed basis sets. 
We will present here a class of invertible maps interrelating the Hopf 
superalgebra $\Uh$ based on quantization of the $r_1$ matrix, and 
its classical analog $\U$. The twist operators {\it vis-{\`a}-vis} 
the above maps will be discussed in the sequel. 
It is shown that a particular map called `minimal twist map'
implements the simplest twist given directly by the factorized form 
of the universal ${\bf R_h}$ matrix of the $\Uh$ algebra. For a 
`non-minimal' map the twist has an additional factor. We evaluate this
twist operator as a series in the deformation parameter.

\setcounter{equation}{0}
\section{A Contraction Process}

A quantization of the classical matrix $r_2$ given in (\ref{eq:clrma})
using a nonlinear basis elements was obtained previously \cite{ACS03}. 
The generating elements $(H, E, F)$ of the $\cUh$ algebra go to the 
classical generators $(h, e, f)$ in the limiting value of the 
deformation parameter: ${\sf h} \rightarrow 0$. Introducing also the 
elements $(X, Y)$ going to $(b_+, b_-)$ respectively in the said 
limit, the commutation relations read
\bea
&&[H, E] = \frac{1}{2}\,(T + T^{-1})\,E , \qquad
[H, F] = - \frac{1}{4}\,(T + T^{-1})\,F 
- \frac{1}{4}\,F\,(T + T^{-1}), \nn\\ 
&&\{E, F\} = - H,\qquad [H, T^{\pm 1}] = T^{\pm 2} - 1,\nn\\
&&[H, Y] = - \frac{1}{2}(T + T^{-1}) Y - \frac{1}{2} Y (T + T^{-1}) 
- \frac{{\sf h}}{4} E (T - T^{-1}) F 
- \frac{\sf h}{4} F (T - T^{-1}) E,\nn\\
&&[T^{\pm 1}, Y] = \pm \frac{{\sf h}}{2} (T^{\pm 1} H + H T^{\pm 1}), 
\qquad E^2 = \frac{T - T^{-1}}{2 \sf h}, \qquad F^2 = - Y,\nn\\
&&[T^{\pm 1}, F] = \pm {\sf h} 
T^{\pm 1} E,\qquad 
[Y, E] = \frac{1}{4} (T + T^{-1}) F + \frac{1}{4} F (T + T^{-1}),
\label{eq:r2alg}
\eea  
where $T^{\pm 1} = \exp({\pm \sf h} X)$. The quantum coproduct $(\Delta)$ 
map for the $\cUh$ algebra reads 
\bea
&&\Delta(H) = H \otimes T^{-1} + T \otimes H 
+ {\sf h} E T^{1/2} \otimes E T^{-1/2},\quad
\Delta(E) = E \otimes T^{-1/2} + T^{1/2} \otimes E,\nn\\    
&&\Delta(F) = F \otimes T^{-1/2} + T^{1/2} \otimes F,
\qquad \Delta(T^{\pm 1}) = T^{\pm 1} \otimes T^{\pm 1},\nn\\
&&\Delta(Y) = Y \otimes T^{-1} + T \otimes Y 
+ \frac{\sf h}{2} E T^{1/2} \otimes T^{-1/2} F
+ \frac{\sf h}{2} T^{1/2} F \otimes E T^{-1/2}.
\label{eq:r2copro}
\eea
The corresponding counit $(\varepsilon)$ and the antipode $(S)$ 
maps are given by
\bea
&&\varepsilon(H) = \varepsilon(E) = \varepsilon(F) 
= \varepsilon(Y) = 0, \qquad \varepsilon(T^{\pm 1}) = 1,\nn\\ 
&&S(H) = - H - {\sf h} E^2, \qquad S(E) = - E, \qquad
S(F) = - F + \frac{\sf h}{2} E,\nn\\
&&S(T^{\pm 1}) = T^{\mp 1}, \qquad 
S(Y) = -Y +\frac{\sf h}{2} H + \frac{{\sf h}^2}{4} E^2.
\label{eq:r2coant}
\eea    
The only primitive element in the above algebra is $X$. We also note 
that the Hopf algebra $\cUh$ has only {\it one} Borel subalgebra 
generated by the elements $(H, E, X)$. It was observed by Kulish
\cite{K98} that the $R_{\sf h}$ matrix in the fundamental 
representation of the quantized algebra $\cUh$ corresponding to the 
classical $r_2$ matrix may be obtained {\it via} a transformation, 
singular in the $q \to 1$ limit, from the corresponding $R_q$ matrix 
in the fundamental representation of the standard $q$-deformed $\Uq$
algebra. The $R_{\sf h}^{1/2;1/2}$ matrix, thus obtained in 
Ref. \cite{K98}, reads    
\beq
R_{\sf h}^{\frac{1}{2};\frac{1}{2}} = 
\left(\begin{array}{ccccccccc}
\noalign{\medskip}1&0&{\sf h}&0&{\sf h}&0
&-{\sf h}&0&\frac {{\sf h}^2}{2}\cr
\noalign{\medskip}0&1&0&0&0&{\sf h}&0&0&0\cr
\noalign{\medskip}0&0&1&0&0&0&0&0&{\sf h}\cr
\noalign{\medskip}0&0&0&1&0&0&0&-{\sf h}&0\cr
\noalign{\medskip}0&0&0&0&1&0&0&0&-{\sf h}\cr
\noalign{\medskip}0&0&0&0&0&1&0&0&0\cr
\noalign{\medskip}0&0&0&0&0&0&1&0&-{\sf h}\cr
\noalign{\medskip}0&0&0&0&0&0&0&1&0\cr
\noalign{\medskip}0&0&0&0&0&0&0&0&1
\label{eq:Rfund}
\end{array}\right).
\eeq

\par

Our task in the present Section is to generalize the above contraction 
procedure for arbitrary representations; and, in particular, we 
construct the $L$ operator corresponding to the Borel subalgebra of 
the Hopf algebra $\cUh$ from the corresponding $L$ operator of the 
standard $q$-deformed $\Uq$ algebra. To this end, we first quote 
the well-known \cite{KR89},\cite{S90} results on the $\Uq$ algebra. 
The quasitriangular quantum Hopf superalgebra $\Uq$, by analogy 
with the classical $\U$ algebra, is generated by three elements 
$({\hat h},\;{\hat e},\;{\hat f})$ obeying the algebraic relations         
\beq
[{\hat h},\;{\hat e}] = {\hat e},\qquad [{\hat h},\;{\hat f}] 
= - {\hat f},\qquad
[{\hat e},\;{\hat f}]=- [h]_{q},
\label{eq:qalg}
\eeq
where $[x]_{q} = ({q^{x}-q^{-x}})/({q-q^{-1}})$. The coalgebraic 
relations are given by  
\bea
&&\Delta\left({\hat h}\right)={\hat h}\otimes 1+1\otimes {\hat h},\qquad
\Delta\left({\hat e}\right)={\hat e}\otimes q^{-{\hat h}/2}
+q^{{\hat h}/2}\otimes {\hat e},\qquad
\Delta\left({\hat f}\right)={\hat f}\otimes q^{-{\hat h}/2}
+q^{{\hat h}/2}\otimes {\hat f},\nn\\
&&\varepsilon({\hat h}) = \varepsilon({\hat e}) 
= \varepsilon({\hat f}) = 0,\qquad 
S({\hat h}) = - {\hat h}, \quad 
S({\hat e}) = - q^{-1/2}\,{\hat e}, \quad 
S({\hat f}) = - q^{1/2}\,{\hat f}. 
\label{eq:qcoalg}
\eea
For convenience, we choose the $(4j+1)$ dimensional irreducible 
representation of the $\Uq$ algebra as follows:
\bea
{\hat h}\;|j\;m> &=& 2 m\;|j\;m>,\qquad 
{\hat e}\;|j\;m> = |j\;m+1/2>,\nn\\
{\hat f}\;|j\;m> &=&- [j + m]_q\,[[j - m + 1/2]]_q\;|j\;m-1/2>\quad
\hbox {for}\;\;j-m\;\;\hbox {integer},\nn\\
\phantom{{\hat f}\;|j\;m> }&=& [[j + m]]_q\,[j - m + 1/2]_q\;|j\;m-1/2>
\quad \hbox {for}\;\;j-m\;\;\hbox {half-integer},
\label{eq:qrep}
\eea
where $[[x]]_q = (q^{x}-(-1)^{2x}q^{-x})/(q^{1/2}+q^{-1/2})$. A related 
numerical quantity used subsequently reads 
$[n]_+ = (-1)^{n-1}\,[[n/2]]_q$, where $n$ is an integer. The 
representation (\ref{eq:qrep}) may be viewed as the $q$-deformation of 
its classical analog (\ref{eq:clrep}). The fundamental representation 
of the algebra (\ref{eq:qalg}) reads
\beq
{\hat h}=\pmatrix{1&0&0\cr 0&0&0\cr 0&0&-1},
\qquad {\hat e}=\pmatrix{0&1&0\cr 0&0&1\cr 0&0&0},
\qquad {\hat f}=\pmatrix{0&0&0\cr -1&0&0\cr 0&1&0}.
\label{eq:funrep}
\eeq
The universal ${\cal R}$-matrix of the $\Uq$ algebra is given by 
\cite{KR89} 
\beq
{\cal R}=q^{{\hat h}\otimes{\hat h}}
\sum_{n=0}^{\infty}\frac{q^{n(n+1)/4}\left(1-q^{-2}\right)^n}
{[n]_+!}\left(q^{{\hat h}/2}{\hat e}\right)^n
\otimes\left(q^{-{\hat h}/2}{\hat f}\right)^n,  
\label{eq:Rq}
\eeq         
where $[n]_+! = [n]_+\,[n-1]_+ \cdots [1]_+, [0]_+! = 1$. 

\par

Following the strategy adopted earlier \cite{ACC98}-\cite{ACC02} for 
constructing the Jordanian deformation of the $sl(N)$ algebra, we
give here a general recipe for obtaining the quantum 
$R_{\sf h}^{j_1;j_2}$ matrix of an arbitrary representation of the  
$\cUh$ algebra. Explicit demonstration is given for the $1/2 \otimes j$
representation, as the relevant $R_{\sf h}^{1/2;j}$ matrix may 
be directly interpreted as the $L$ operator corresponding to 
the Borel subalgebra of the $\cUh$ algebra. But our construction may be 
obviously generalized. The primary ingredient for our method is the 
$R_q^{1/2;j}$ matrix of the $\Uq$ algebra in the $1/2 \otimes j$ 
representation. A suitable similarity transformation is performed on 
this $R_q^{1/2;j}$ matrix. The transforming matrix is singular in the 
$q \rightarrow 1$ limit. {\it For the transformed matrix, the 
singularities, however, systematically cancel yielding a well-defined 
construction}. The transformed matrix, directly furnishes the 
$R_{\sf h}^{1/2;j}$ for the nonstandard $\cUh$ algebra. Interpreting, 
as mentioned above, the $R_{\sf h}^{1/2;j}$ obtained here as the $L$ 
operator corresponding to the Borel subalgebra of the $\cUh$ algebra, 
we use the standard FRT procedure \cite{FRT90} to reconstruct the full 
Hopf structure of the said Borel subalgebra presented 
in (\ref{eq:r2alg})-(\ref{eq:r2coant}). The $R_q^{1/2;j}$  matrix 
of the tensored $1/2 \otimes j$ representation of the $\Uq$ 
algebra reads
\beq        
R_q^{\frac{1}{2};j} = \left(\begin{array}{ccc}
q^{\hat h}&-\omega q^{{\hat h}/2}{\hat f}&
- \omega \left(1+q^{-1}\right){\hat f}^2\\
\noalign{\medskip}0&1&\omega q^{-{(\hat h+1)}/2}{\hat f}\\ 
\noalign{\medskip}0&0&q^{-{\hat h}} \end{array}\right),
\label{eq:Rqfunarb}
\eeq
where $\omega = q - q^{-1}$. We now introduce a transforming matrix 
$M$, singular in the $q \rightarrow 1$ limit, as
\beq 
M = E_{q^2}(\eta {\hat e}^2),
\label{eq:singop}
\eeq
where
\beq 
E_{q^2}(x) = \sum_{n=0}^{\infty} \frac{x^n}{[n]_{q^2}!},\qquad 
\eta = \frac{\sf h}{q^2 -1}.
\label{eq:Eexp}
\eeq 
For any finite value of $j$ the series (\ref{eq:Eexp}) may be terminated 
after setting ${\hat e}^{4j+1} = 0$; but, we proceed quite generally.  
The $R_q^{j_1;j_2}$ matrix of the $\Uq$ algebra may now be subjected 
to a similarity transformation followed by a limiting process
\beq
\tilde{R}_{\sf h}^{j_1;j_2} \equiv \lim_{q \to 1} 
\left[\left( M_{j_1}^{-1} \otimes M_{j_2}^{-1} \right)\,
R_q^{j_1;j_2}\,\left( M_{j_1} \otimes M_{j_2} \right)\right].
\label{eq:gencon}   
\eeq
In the followings we will present explicit results for the operator 
$\tilde{R}_{\sf h}^{1/2;j}$. In performing the similarity 
transformation (\ref{eq:gencon}) we may choose any suitable operator 
ordering. Specifically, starting from left we maintain the order 
${\hat e} \prec {\hat h} \prec {\hat f}$. In our calculation a class of 
operators 
\beq
{\cal T}_{(\alpha)} = 
(E_{q^2}(\eta {\hat e}^2))^{-1}\,E_{q^2}(q^{2 \alpha}\eta {\hat e}^2)
\label{eq:Tauop}
\eeq
satisfying
\beq
(E_{q^2}(\eta {\hat e}^2))^{-1}\,q^{\alpha {\hat h}}\,
(E_{q^2}(\eta {\hat e}^2))
= {\cal T}_{(\alpha)} q^{\alpha {\hat h}},\qquad 
{\cal T}_{(\alpha + \beta)} q^{(\alpha + \beta) {\hat h}} 
= {\cal T}_{(\alpha)} q^{\alpha {\hat h}}\;
{\cal T}_{(\beta)} q^{\beta {\hat h}}, 
\label{eq:Trule}
\eeq
play an important role. To evaluate $q \to 1$ limiting value of the 
operator ${\cal T}_{(\alpha)}$, we use the identity
\beq
E_{q^2}(q^2\,\eta\,{\hat e}^2) - E_{q^2}(q^{-2}\,\eta\,{\hat e}^2)
= \eta\;(q^2 - q^{-2})\;{\hat e}^2\;E_{q^2}(\eta {\hat e}^2),
\label{eq:Eiden}
\eeq
which follows from the series ({\ref{eq:Eexp}). The identity 
(\ref{eq:Eiden}) may be rephrased as
\beq
{\cal T}_{(1)} - {\cal T}_{(-1)} = \eta\;(q^2 -q^{-2})\;{\hat e}^2.
\label{eq:Tiden}
\eeq     
Evaluating term by term, the limiting values of 
${\cal T}_{(\pm 1)}|_{q \rightarrow 1}\;
\left(\equiv {\tilde T}_{(\pm 1)}\right)$ are found to be {\em finite};
and, for these finite operators the second equation in 
(\ref{eq:Trule}) suggests that 
\beq
{\tilde T}_{(\pm \alpha)} = \left({\tilde T}_{(\pm 1)}\right)^{\alpha},
\label{eq:Talpha}     
\eeq
where ${\tilde T}_{(\alpha)} = \lim_{q \to 1}{\cal T}_{(\alpha)}$. 
Writing ${\tilde T}_{(\pm 1)} = {\tilde T}^{\pm 1}$ henceforth, we 
immediately observe that in the $q \to 1$ limit, the identity 
(\ref{eq:Tiden}) assumes to the form 
\beq
{\tilde T} - {\tilde T}^{-1} = 2{\sf h} e^2,
\label{eq:Tq1iden}
\eeq
which may be solved as 
\beq
{\tilde T}^{\pm 1} = \pm {\sf h} e^2 + \sqrt{1 + {\sf h}^2 e^4}. 
\label{eq:Tsolve}
\eeq
This is our crucial result. Two other operator identities playing 
key roles are listed below:
\bea
{\hat f}\,{\hat e}^{2n} &=& {\hat e}^{2n}\,{\hat f} 
- \frac{q}{q + 1}\,\{n\}_{q^2}\,{\hat e}^{2n-1}\,{\hat t}
- \frac{1}{q + 1}\,\{n\}_{q^{-2}}\,{\hat e}^{2n-1}\,{\hat t}^{-1},\\
{\hat f}^2\,{\hat e}^{2n} &=& {\hat e}^{2n}\,{\hat f}^2
+ q \frac{q-1}{q+1}\,\{n\}_{q^2}\,{\hat e}^{2n-1}\,{\hat t}\,{\hat f} 
- q^{-1} \frac{q-1}{q+1}\,\{n\}_{q^{-2}}\,
{\hat e}^{2n-1}\,{\hat t}^{-1}\,{\hat f}\nn\\ 
&\phantom{=}& + \frac{q}{q+1}\,\left(\frac{1}{\omega}\{n\}_{q^4} 
- q^2\,\frac{q-1}{q+1}\,\frac{\{n-1\}_{q^2}\,\{n\}_{q^2}}{\{2\}_{q^2}}
\right)\,{\hat e}^{2(n-1)}\,{\hat t}^2\nn\\
&\phantom{=}& - \frac{1}{q+1}\,\left(\frac{1}{\omega}\{n\}_{q^{-4}} 
- q^{-2}\,\frac{q-1}{q+1}\,
\frac{\{n-1\}_{q^{-2}}\,\{n\}_{q^{-2}}}{\{2\}_{q^{-2}}}\right)
\,{\hat e}^{2(n-1)}\,{\hat t}^{-2}\nn\\
&\phantom{=}& - \frac{q}{(q+1)^3}\,\left(q\{n\}_{q^2} 
+ \{n\}_{q^{-2}}\right)\,{\hat e}^{2(n-1)},
\label{eq:OPE}
\eea
where $\{x\}_q = (1 - q^x)/(1 - q)$ and 
${\hat t}^{\pm 1} = q^{\pm {\hat h}}$. Using the above two identities 
systematically and passing to the limit $q \to 1$, it may be shown that 
in our construction of the operator ${\tilde R}_{\sf h}^{\frac{1}{2};j}$    
{\it via} (\ref{eq:gencon}), {\em all singularities cancel} yielding a 
well-defined answer 
\beq
{\tilde R}_{\sf h}^{\frac{1}{2};j} =\left(\begin{array}{ccc}
{\tilde T}&{\sf h} {\tilde T}^{\frac{1}{2}} e&- {\sf h} {\tilde H}
+\frac{\sf h}{4}\left({\tilde T} - {\tilde T}^{-1}\right)\cr 
\noalign{\medskip}0&1&-{\sf h} {\tilde T}^{- \frac{1}{2}} e\cr 
\noalign{\medskip}0&0&{\tilde T}^{-1}
\end{array}\right),
\label{eq:contR}
\eeq
where ${\tilde H} = \frac{1}{2} ({\tilde T} + {\tilde T}^{-1})\;h 
= \sqrt{1 + {\sf h}^2 e^4}\;h$. One way of interpreting 
(\ref{eq:contR}) is to consider it a recipe for obtaining the 
finite dimensional $R_{\sf h}$ matrices of the $\cUh$ algebra.
For instance, using the classical $j = 1$ representation given in 
(\ref{eq:clrep}) we obtain the $R_{\sf h}^{1/2;1} \left(= 
{\tilde R}_{\sf h}^{1/2;1}\right)$ as follows:
\beq
R_{\sf h}^{\frac{1}{2};1} = \left(\begin {array}{ccccccccccccccc}
1&0&{\sf h}&0&\frac{{\sf h}^2}2&0&{\sf h}&0&\frac{{\sf h}^2}2&
0&-2{\sf h}&0&\frac{{\sf h}^2}2&0&{\sf h}^3\\
\noalign{\medskip}0&1&0&{\sf h}&0&0&0&{\sf h}&0&\frac{{\sf h}^2}2&0&
-{\sf h}&0&\frac{{\sf h}^2}2&0\\
\noalign{\medskip}0&0&1&0&{\sf h}&0&0&0&{\sf h}&0&0&0&0&0&
\frac{{\sf h}^2}2\\
\noalign{\medskip}0&0&0&1&0&0&0&0&0&{\sf h}&0&0&0&{\sf h}&0\\
\noalign{\medskip}0&0&0&0&1&0&0&0&0&0&0&0&0&0&2{\sf h}\\
\noalign{\medskip}0&0&0&0&0&1&0&0&0&0&0&-{\sf h}&0&
\frac{{\sf h}^2}2&0\\
\noalign{\medskip}0&0&0&0&0&0&1&0&0&0&0&0&-{\sf h}&0&
\frac{{\sf h}^2}2\\
\noalign{\medskip}0&0&0&0&0&0&0&1&0&0&0&0&0&-{\sf h}&0\\
\noalign{\medskip}0&0&0&0&0&0&0&0&1&0&0&0&0&0&-{\sf h}\\
\noalign{\medskip}0&0&0&0&0&0&0&0&0&1&0&0&0&0&0\\
\noalign{\medskip}0&0&0&0&0&0&0&0&0&0&1&0&-{\sf h}&0&
\frac{{\sf h}^2}2\\
\noalign{\medskip}0&0&0&0&0&0&0&0&0&0&0&1&0&-{\sf h}&0\\
\noalign{\medskip}0&0&0&0&0&0&0&0&0&0&0&0&1&0&-{\sf h}\\
\noalign{\medskip}0&0&0&0&0&0&0&0&0&0&0&0&0&1&0\\
\noalign{\medskip}0&0&0&0&0&0&0&0&0&0&0&0&0&0&1\end{array}\right).
\label{eq:Rhhalfone}
\eeq
The matrix (\ref{eq:contR}) may also be interpreted as the $L$ 
operator of the $\cUh$ algebra. To this end, we first use the 
following invertible map of the quantum $\cUh$ algebra 
(\ref{eq:r2alg}) on the classical algebra (\ref{eq:clalg}): 
\beq
E = e,\quad H = {\tilde H},\quad 
F = f + \frac{\sf h}{4}\,\left(\frac{{\tilde T} - 1}
{{\tilde T} + 1}\right)\,e 
- \frac{\sf h}{2}\,\left(\frac{{\tilde T} - 1}
{{\tilde T} + 1}\right)\,eh,\quad
T = {\tilde T},\quad Y = F^2. 
\label{eq:contrmap} 
\eeq
The map (\ref{eq:contrmap}) satisfies the algebraic relations 
(\ref{eq:r2alg}). The general structure of maps relating the 
$\cUh$ algebra, obtained by quantizing the classical $r_2$ matrix,
and the classical $\U$ algebra has been discussed before 
\cite{ACS03}. Using the map (\ref{eq:contrmap}) the operator 
(\ref{eq:contR}) may be recast as
\beq
L \equiv R_{\sf h}^{\frac{1}{2};j} = \left(\begin{array}{ccc}
T&{\sf h} T^{\frac{1}{2}} E&- {\sf h} H
+\frac{\sf h}{4} (T - T^{-1})\cr 
\noalign{\medskip}0&1&-{\sf h} T^{- \frac{1}{2}} E \cr 
\noalign{\medskip}0&0&T^{-1}
\end{array}\right).
\label{eq:Lhalf}
\eeq
The above $L$ operator allows immediate construction of the full Hopf 
structure of the Borel subalgebra of the $\cUh$ algebra {\it via} the 
standard FRT formalism \cite{FRT90}. The algebraic relations for the 
generators ($E, T^{\pm 1}, H$) of the Borel subalgebra is given by 
\beq
R_{\sf h}^{\frac{1}{2};\frac{1}{2}}\,L_{1}\,L_{2}  
= L_{2}\,L_{1}\,R_{\sf h}^{\frac{1}{2};\frac{1}{2}},
\label{eq:RLL}
\eeq  
where $\ZZ_{2}$ graded tensor product has been used in defining the
operators: $L_{1} = L \otimes {\mathbf I},\, 
L_{2} = {\mathbf I} \otimes L$. The coalgebraic properties of the 
said Borel subalgebra may be succinctly expressed as
\beq
\Delta(L) = L \dot{\otimes} L,\qquad
\varepsilon(L) = {\mathbf I},\qquad S(L) = L^{-1},
\label{eq:Lcoalg} 
\eeq 
where $L^{-1}$ is given by
\beq
L^{-1} = \left(\begin{array}{ccc}
T^{-1}& -{\sf h} T^{- \frac{1}{2}} E & 
{\sf h} H + \frac{\sf h}{4}\left(T-T^{-1}\right)\\ 
\noalign{\medskip}0&1&{\sf h} T^{\frac{1}{2}} E\\
\noalign{\medskip}0&0&T
\end{array}\right).
\label{eq:Linv}
\eeq
This completes our construction of the Hopf structure of the Borel 
subalgebra of the $\cUh$ algebra, obtained by deforming the $r_2$
matrix, by employing the contraction scheme described earlier.  
Our recipe (\ref{eq:gencon}) for obtaining  
$R_{\sf h}^{j_1;j_2}$ matrices for $j_1 \otimes j_2$ 
representation of the $\cUh$ algebra may be continued arbitrarily. 
The matrices such as $R_{\sf h}^{1;j}$ may be interpreted
as higher dimensional $L$ operators \cite{JJ96} obeying duality
relations with higher representations of ${\cal T}$ matrices.  

\setcounter{equation}{0}
\section{A nonlinear realization of the quantized $\Uh$ algebra 
corresponding to the $r_1$ matrix}

The classical $r_1$ matrix has been quantized earlier \cite{CK98}
using the inclusion $sl(2) \subset osp(2|1)$. These authors
have expressed the resultant triangular deformed $osp(2|1)$ algebra
in terms of the classical basis set. On the other hand Ohn \cite{O92}
employed a nonlinear basis set to formulate the Jordanian deformed
${\cal U}_{\sf h}(sl(2))$ algebra. Consequently, Ohn's   
${\cal U}_{\sf h}(sl(2))$ algebra \cite{O92} do not directly appear
as a Hopf subalgebra of the deformed $osp(2|1)$ considered
in Ref. \cite{CK98}. Moreover, if the algebraic relations are described 
in terms of the undeformed classical basis set, the coproduct structures
tend to be complicated in nature.   

In our work, we present the quantized Hopf structure corresponding to 
the classical $r_1$ matrix in terms of nonlinear basis elements. Our
algebra {\it explicitly includes} Ohn's ${\cal U}_{\sf h} (sl(2))$
algebra as a Hopf subalgebra. The coproduct structure we obtain
is {\it considerably simple}. To distinguish the deformed $osp(2|1)$ 
algebra, its generators and the deformation parameter considered in 
the present Section from the corresponding objects presented in 
Section 2, we express them in boldfaced notations. The quantized Hopf 
algebra $\Uh$ corresponding to the classical $r_1$ matrix is generated 
by the elements $({\bf H, E, F, X, Y})$. Their classical analogs are 
$(h, e, f, b_{\pm})$ respectively. The elements ${\bf T}^{\pm 1} = 
\exp ({\pm \bf hX})$ are also introduced. The deformation parameter is 
denoted by ${\bf h}$. The Hopf structure of the $\Uh$ algebra is 
obtained by maintaining the following properties: ({\sf i}) In the 
classical limit the quantum coproduct map conforms to the classical
cocommutator. ({\sf ii}) The coproduct map is a homomorphism of 
the algebra, and it satisfies the coassociativity constraint. 
({\sf iii}) Generator ${\bf X}$ is the only primitive element.          
The commutation relations of the $\Uh$ algebra reads 
\bea
&& [{\bf H},\,{\bf E}] = \frac {1}{2}
({\bf T}+{\bf T}^{-1}){\bf E},\nn\\
&& [{\bf H},\,{\bf F}] = 
-\frac {1}{4} ({\bf T}+{\bf T}^{-1}){\bf F}
-\frac {1}{4}{\bf F} ({\bf T}+{\bf T}^{-1})
-\frac{\bf h}{8} (({\bf T}-{\bf T}^{-1}) {\bf H}
+ {\bf H} ({\bf T}-{\bf T}^{-1})) {\bf E} \nn\\
&&\phantom{[{\bf H},\,{\bf F}]=} - \frac{\bf h}{8} {\bf E} 
(({\bf T}-{\bf T}^{-1}) {\bf H}
+ {\bf H}({\bf T}-{\bf T}^{-1})),\nn\\
&& \{{\bf E},\,{\bf F}\} = 
-\frac {1}{4} ({\bf T}+{\bf T}^{-1}) {\bf H} 
- \frac {1}{4} {\bf H} ({\bf T}+{\bf T}^{-1}),\qquad
[{\bf H},\,{\bf T}^{\pm 1}] = {\bf T}^{\pm 2} - 1,\nn\\
&&[{\bf H},\,{\bf Y}] = 
- \frac{1}{2} ({\bf T} + {\bf T}^{-1}) {\bf Y}
- \frac{1}{2} {\bf Y} ({\bf T} + {\bf T}^{-1}),\qquad
[{\bf T}^{\pm 1},\,{\bf Y}] = \pm \frac{\bf h}{2} 
({\bf T}^{\pm 1} {\bf H} + {\bf H } {\bf T}^{\pm 1}),\nn\\ 
&&{\bf E}^2=\frac {1}{2{\bf h}} ({\bf T}-{\bf T}^{-1}),
\qquad [{\bf Y}, {\bf E}] = {\bf F},\qquad 
[{\bf T}^{\pm 1}, {\bf F}] = 
\pm \frac{\bf h}{2} ({\bf T}^{\pm 2} + 1) {\bf E},\nn\\
&& {\bf F}^2 = -{\bf Y}+\frac{\bf h}{8}
({\bf T}-{\bf T}^{-1}) {\bf H}^2
+ \frac{\bf h}{4} ({\bf T}-{\bf T}^{-1}) {\bf E F}
+ \frac{3{\bf h}}{16}({\bf T}^2-{\bf T}^{-2}) {\bf H}
+ \frac{\bf h}{4} ({\bf T}-{\bf T}^{-1})\nn\\
&&\phantom{{\bf F}^2 =} + \frac{9{\bf h}}{128}
({\bf T}-{\bf T}^{-1})^3,\nn\\
&&[{\bf F},\;{\bf Y}] = 
\frac{\bf h}{4} ({\bf T}-{\bf T}^{-1}) {\bf F} 
+ \frac{\bf h}{2} ({\bf T}-{\bf T}^{-1}) {\bf E Y}  
- \frac{{\bf h}^2}{4}{\bf E H}^2
- \frac{3{\bf h}^2}{8} ({\bf T}+{\bf T}^{-1}) {\bf E H} 
- \frac{{\bf h}^2}{2} {\bf E}\nn\\
&&\phantom{[{\bf F},\;{\bf Y}] =} - \frac{15{\bf h}^2}{64}
({\bf T}-{\bf T}^{-1})^2 {\bf E}
\label{eq:r1alg}
\eea
and the corresponding coalgebraic structure is given by
\bea
&&\Delta({\bf H}) = {\bf H} \otimes {\bf T} 
+ {\bf T}^{-1} \otimes {\bf H},\qquad\qquad\qquad\qquad 
\Delta({\bf E}) = {\bf E} \otimes {\bf T}^{-1/2} 
+ {\bf T}^{1/2}\otimes {\bf E},\nn\\
&&\Delta ({\bf F}) = {\bf F} \otimes {\bf T}^{1/2}
+ {\bf T}^{-1/2}\otimes {\bf F}+ \frac{\bf h}{4} {\bf T}^{-1} {\bf E} 
\otimes \left({\bf T}^{-1/2}{\bf H} + {\bf H T}^{-1/2}\right)
- \frac{\bf h}{4} \left ({\bf T}^{1/2} {\bf H}
+ {\bf H T}^{1/2}\right)\otimes {\bf T E},\nn\\ 
&&\Delta({\bf T}^{\pm 1}) = {\bf T}^{\pm 1} \otimes {\bf T}^{\pm 1},
\qquad\qquad\qquad\qquad \Delta({\bf Y}) = {\bf Y} \otimes {\bf T} 
+ {\bf T}^{-1} \otimes {\bf Y},\nn\\
&&\varepsilon({\bf H}) = \varepsilon({\bf E}) = \varepsilon({\bf F}) 
= \varepsilon({\bf Y}) = 0,\qquad\qquad\qquad 
\varepsilon({\bf T}^{\pm 1}) = 1,\nn\\
&&S({\bf H}) = - {\bf H} + 2 {\bf h} {\bf E}^2,\qquad\qquad 
S({\bf E}) = - {\bf E},\qquad\qquad
S({\bf F}) = - {\bf F} 
- \frac{\bf h}{2}\left({\bf T} + {\bf T}^{-1}\right) {\bf E},\nn\\
&&S({\bf T}^{\pm 1}) = {\bf T}^{\mp 1},\qquad\qquad\qquad\qquad
S({\bf Y}) = - {\bf Y} - {\bf h H} + {\bf h}^2 {\bf E}^2.
\label{eq:r1coalg}
\eea
All the Hopf superalgebra axioms can be verified by direct calculation. 
The universal ${\bf R}_{\bf h}$ matrix of $\Uh$ is of the factorized 
form \cite{BH96}:  
\beq
{\bf R_h} = G_{21}^{-1} G,\qquad G = \exp\left({\bf h}\,{\bf T H} 
\otimes {\bf X}\right),  
\label{eq:r1univR}
\eeq
which coincides with the universal ${\cal R}_{\sf h}$ matrix 
of the ${\cal U}_{\sf h}(sl(2))$ subalgebra \cite{BH96} involving the 
highest weight root vector.

\par

Before discussing the general structure of a class of invertible maps 
of the $\Uh$ algebra on the classical $\U$ algebra, we notice the
comultiplication map of a set of three operators 
\beq
{\bf T}^{-1/2} {\bf E},\quad {\bf T H}, 
\quad {\bf T}^{1/2} {\bf F} 
+ \frac{\bf h}{8} {\bf T}^{1/2} ({\bf T} - {\bf T}^{-1}) {\bf E}
- \frac{\bf h}{2} {\bf T}^{1/2} {\bf E H},
\label{eq:mtmop}
\eeq
when acted by the twist operator corresponding to the factorized 
form of the universal ${\bf R_h}$ matrix, reduce to the classical 
cocommutative coproduct:
\beq
G \Delta({\cal X}) G^{-1} = {\cal X} \otimes \mathbf{I} 
+ \mathbf{I} \otimes {\cal X},  
\label{eq:mtmtwist}
\eeq 
where ${\cal X}$ is an element of the set (\ref{eq:mtmop}). From the 
commutation rules (\ref{eq:r1alg}) it becomes evident that the 
operators in the set (\ref{eq:mtmop}) satisfy the classical algebra
generated by $(e, h, f)$ respectively. A general discussion of the 
invertible maps between the $\cUh$ algebra based on quantization of the 
$r_2$ matrix, and the classical $\U$ algebra was given in 
Ref. \cite{ACS03}. Our present discussion of maps interrelating the
$\Uh$ algebra, obtained {\it via} quantization of the $r_1$ matrix,
and its classical analog $\U$ follow the same pattern. The details of 
the construction are, however, quite different. A short description
of the maps in the present case is given below.

\par  
 
The quantum $\Uh$ algebra may be mapped on the classical $\U$ algebra
by using a general ansatz as follows:
\bea
{\bf E} &=& \varphi_1 (b_{+}) e,\qquad 
{\bf H} = \varphi_2(b_{+}) h,\nn\\
{\bf F} &=& \varphi_3(b_{+}) f + u_1(b_{+}) e + u_2(b_{+}) e h,
\label{eq:ansatz}
\eea
where the 'mapping functions' ($\varphi_1, \varphi_2, \varphi_3; 
u_1, u_2$) depend only on the classical generator $b_{+}$. In the 
classical limit ${\bf h} \to 0$ the above functions satisfy the property: 
($\varphi_1, \varphi_2, \varphi_3; u_1, u_2$) $\to$ ($1, 1, 1; 0, 0$). 
The operators ${\bf T}^{\pm 1}$ may now be expressed as 
\beq
{\bf T}^{\pm 1} = \pm {\bf h} b_{+} (\varphi_1(b_{+}))^2
+ \sqrt{ 1 + {\bf h}^2 b_{+}^2 (\varphi_1(b_{+}))^4}.
\label{eq:TiTin}
\eeq 
Substituting the ansatz (\ref{eq:ansatz}) in the defining relations
(\ref{eq:r1alg}) for the $\Uh$ algebra we, for a {\it given} function 
$\varphi_1$, obtain a set of {\it six} nonlinear equations for 
{\it four} unknown functions:
\bea
&&\left(\varphi_1(b_{+}) + 2 b_{+} \varphi_1^{\prime}(b_{+})\right) 
\varphi_2(b_{+}) 
- \sqrt{ 1 + {\bf h}^2 b_{+}^2 (\varphi_1(b_{+}))^4}\,
\varphi_1(b_{+}) = 0,\nn\\
&&2 b_{+} \varphi_2(b_{+}) \varphi_3^{\prime}(b_{+})
- \varphi_2(b_{+}) \varphi_3(b_{+}) 
+ \sqrt{ 1 + {\bf h}^2 b_{+}^2 (\varphi_1(b_{+}))^4}\, 
\varphi_3(b_{+})= 0,\nn\\
&&2 b_{+} \varphi_2(b_{+}) u_1^{\prime}(b_{+})
+ \left(\varphi_2(b_{+}) + \sqrt{1 
+ {\bf h}^2 b_{+}^2 (\varphi_2(b_{+}))^4}\right)\,u_1(b{+})\nn\\
&&\qquad\qquad + {\bf h}^2 b_{+} 
\sqrt {1 + {\bf h}^2 b_{+}^2 (\varphi_1(b_{+}))^4}\,
(\varphi_1(b_{+}))^3 = 0,\nn\\ 
&&\left(\varphi_2(b_{+}) - 2 b_{+}\varphi_2^{\prime}(b_{+})
+ \sqrt {1 + {\bf h}^2 b_{+}^2 (\varphi_1(b_{+}))^4} \right) u_2(b_{+}) 
+ \varphi_2^{\prime}(b_{+}) \varphi_3(b_{+})\nn\\ 
&&\qquad\qquad
+ 2 b_{+} \varphi_2(b_{+}) u_2^{\prime}(b_{+}) 
+ {\bf h}^2 b_{+} (\varphi_1(b_{+}))^3 \varphi_2(b_{+}) = 0,\nn\\
&&\varphi_1(b_{+}) \left(2 b_{+} u_2(b_{+}) - \varphi_3(b_{+})\right) 
+ \sqrt {1 + {\bf h}^2 b_{+}^2 (\varphi_1(b_{+}))^4}\, 
\varphi_2(b_{+}) = 0,\nn\\ 
&&b_{+} \varphi_1(b_{+}) \left(2 u_1(b_{+}) + u_2(b_{+})\right)
- b_{+} \varphi_1^{\prime}(b_{+}) \left( \varphi_3(b_{+}) 
- 2 b_{+} u_2(b_{+})\right) 
+ {\bf h}^2 b_{+}^2 (\varphi_1(b_{+}))^4 = 0. 
\label{eq:fordiff}
\eea
Maintaining the classical limit the above set of equations may be
consistently solved as follows:
\bea
&&\varphi_2(b_{+}) = 
\frac{\sqrt{ 1 + {\bf h}^2 b_{+}^2 (\varphi_1(b_{+}))^4}\,
\varphi_1(b_{+})}
{\varphi_1(b_{+}) + 2 b_{+} \varphi_1^{\prime}(b_{+})},\qquad
\varphi_3(b_{+}) = \frac{1}{\varphi_1(b_{+})},\nn\\ 
&&u_1(b_{+}) = - \frac{{\bf h}^2}{4}\,b_{+} (\varphi_1(b_{+}))^3,\qquad 
u_2(b_{+}) = 
\frac{1 - \sqrt{ 1 + {\bf h}^2 b_{+}^2 (\varphi_1(b_{+}))^4}\,
\varphi_2(b_{+})}{2 b_{+} \varphi_1(b_{+})}.
\label{eq:forsol}
\eea
Inverse maps expressing the classical generators in terms of the 
relevent quantum generators are obtained by assuming the anstaz
\bea
e &=& \psi_1({\bf T}),\qquad\;\;h = \psi_2({\bf T})\,{\bf H},\nn\\
f &=& \psi_3({\bf T})\,{\bf F} + w_1({\bf T})\,{\bf E} 
+ w_2({\bf T})\,{\bf E H},
\label{eq:invan}
\eea  
where $(\psi_1, \psi_2, \psi_3; w_1, w_2)$ are functions of ${\bf T}$ 
obeying the limiting property: $(\psi_1, \psi_2, \psi_3; w_1, w_2) 
\to (1, 1, 1; 0, 0)$ as ${\bf h} \to 0$. The differential equations
obeyed by the 'mapping functions' introduced in (\ref{eq:invan}) 
read
\bea
&&2 ({\bf T}^2 - 1)\,\psi_1^{\prime}({\bf T}) \psi_2({\bf T})     
+ ({\bf T} + {\bf T}^{-1})\,\psi_1({\bf T}) \psi_2({\bf T})     
- 2 \psi_1({\bf T}) = 0,\nn\\     
&&2 ({\bf T}^2 - 1)\,\psi_2({\bf T}) \psi_3^{\prime}({\bf T})     
- ({\bf T} + {\bf T}^{-1})\,\psi_2({\bf T}) \psi_3({\bf T})     
+2 \psi_3({\bf T}) = 0,\nn\\     
&&\psi_2({\bf T})\,\left(2 ({\bf T} + {\bf T}^{-1})\,w_1({\bf T})     
+ 4 ({\bf T}^2 - 1)\,w_1^{\prime}({\bf T})     
- {\bf h} ({\bf T}^2 - {\bf T}^{-2})\,\psi_3({\bf T})\right)
+ 4 w_1({\bf T}) = 0,\nn\\
&&({\bf T} - {\bf T}^{-1})\,\psi_2({\bf T})\,
\left(2{\bf T}\,w_2^{\prime}({\bf T}) - {\bf h}\,\psi_3({\bf T})\right)     
+ {\bf h} ({\bf T}^2 + 1)\,\psi_2^{\prime}({\bf T})\,
\psi_3({\bf T})\nn\\     
&&\qquad\qquad 
+ \left(({\bf T} + {\bf T}^{-1})\,\psi_2({\bf T}) - 2 ({\bf T}^2 - 1)
\,\psi_2^{\prime}({\bf T}) + 2\right) w_2({\bf T}) = 0,\nn\\     
&&2 ({\bf T} - {\bf T}^{-1})\,\psi_1({\bf T}) w_2({\bf T})     
- {\bf h} ({\bf T} + {\bf T}^{-1})\,\psi_1({\bf T})\,\psi_3({\bf T})     
+ 2 {\bf h} \psi_2({\bf T}) = 0,\nn\\
&&({\bf T} - {\bf T}^{-1})\,\psi_1({\bf T})\,\left(4 w_1({\bf T}) 
+  ({\bf T} + {\bf T}^{-1})\,w_2({\bf T}) 
- {\bf h} ({\bf T} - {\bf T}^{-1})\,\psi_3({\bf T})\right)\nn\\ 
&&\qquad\qquad      
+ {\bf T} ({\bf T} - {\bf T}^{-1})\,\psi_1^{\prime}({\bf T})  
\left(2 ({\bf T} - {\bf T}^{-1})\,w_2({\bf T}) 
- {\bf h} ({\bf T} + {\bf T}^{-1})\,\psi_3({\bf T})\right) = 0. 
\label{eq:backdiff}
\eea     
Treating the function $\psi_1({\bf T})$ as known, and maintaining the 
limiting properties the remaining functions may be solved uniquely:
\bea
\psi_2({\bf T}) &=& \frac{2 \psi_1({\bf T})}     
{({\bf T} + {\bf T}^{-1})\,\psi_1({\bf T}) 
+ 2 ({\bf T}^2 - 1)\,\psi_1^{\prime}({\bf T})},\qquad\qquad     
\psi_3({\bf T}) = \frac{1}{\psi_1({\bf T})},\nn\\
w_1({\bf T}) &=& \frac{{\bf h} ({\bf T} - {\bf T}^{-1})}
{8 \psi_1({\bf T})},\qquad\qquad\quad
w_2({\bf T}) = {\bf h}\,\frac{{\bf T} + {\bf T}^{-1} 
- 2 \psi_2({\bf T})}{2 ({\bf T} - {\bf T}^{-1})\,\psi_1({\bf T})}. 
\label{eq:backsol}
\eea

\par

The general structure of the twisting elements corresponding to the 
given maps may be described as follows. Let $m$ be a deformation map
and $m^{-1}$ be its inverse:
\beq
m: ({\bf E, H, F}) \to (e, h, f),\qquad m^{-1}: (e, h, f)
\to ({\bf E, H, F}).
\label{eq:mminv}
\eeq
The classical $(\Delta_0)$ cocommutative and the quantum $(\Delta)$
non-cocommutative coproducts are related \cite{D90} by the twisting 
element as 
\beq
{\cal G}\,\Delta \circ m^{-1} (\phi)\,{\cal G}^{-1} 
= (m^{-1} \otimes m^{-1}) \circ \Delta_0(\phi)\qquad
\forall \phi \in \U,
\label{eq:Drinfeld}   
\eeq
where the twisting element $ {\cal G} \in \Uh^{\otimes 2}$  
satisfies the cocycle condition 
\beq
({\cal G} \otimes {\bf I})\,((\Delta \otimes \hbox {id}) {\cal G}) 
= ({\bf I} \otimes {\cal G})\,((\hbox {id} \otimes \Delta) {\cal G}).
\label{eq:cocycle}
\eeq 
Similarly the classical $(S_0)$ and the quantum $(S)$ antipode maps   
are related as follows:
\beq
{\bf g}\,S \circ m^{-1} (\phi)\,{\bf g}^{-1} 
= m^{-1} \circ S_0(\phi),\qquad
{\bf g} \in \Uh.
\label{eq:Srel}
\eeq
The transforming operator ${\sf g}$ for the antipode map may be 
expressed in terms of the twist operator ${\cal G}$ as
\beq
{\bf g} = \mu \circ (\hbox{id} \otimes S) {\cal G},
\label{eq:antwi}
\eeq 
 where $\mu$ is the multiplication map.

\par

The first map considered here plays a key role in the present 
construction. With the choice      
\beq
\varphi_1(b_{+}) = (1 - 2 {\bf h} b_{+})^{-1/4},\qquad
\psi_1({\bf T}) = {\bf T}^{-1/2},
\label{eq:mtmchoice}
\eeq     
we obtain the following direct map
\bea
{\bf E} &=& (1 - 2 {\bf h} b_{+})^{-1/4}\,e,\quad
{\bf H} = \sqrt{(1 - 2 {\bf h} b_{+})}\,h,\quad
{\bf T}^{\pm 1} = (1 - 2 {\bf h} b_{+})^{\mp 1/2},\nn\\
{\bf F} &=& (1 - 2 {\bf h} b_{+})^{1/4}\,f - \frac{{\bf h}^2}{4} b_{+} 
(1 - 2 {\bf h} b_{+})^{-3/4}\,e
+ \frac{\bf h}{2}\,(1 - 2 {\bf h} b_{+})^{1/4}\,eh
\label{eq:mtmdir}
\eea
and its inverse
\beq
e = {\bf T}^{-1/2} {\bf E},\quad h = {\bf T H}, 
\quad f = {\bf T}^{1/2} {\bf F} 
+ \frac{\bf h}{8} {\bf T}^{1/2} ({\bf T} - {\bf T}^{-1}) {\bf E}
- \frac{\bf h}{2} {\bf T}^{1/2} {\bf E H}.
\label{eq:mtminv}
\eeq
It may be observed from (\ref{eq:mtmtwist}) that the operator $G$
corresponding to the factorized form of the universal ${\bf R_h}$
matrix given in (\ref{eq:r1univR}) plays the role of the twist operator
${\cal G}$ for the map (\ref{eq:mtmdir}) and its inverse. In this sense 
we refer to it as the 'minimal twist map'. The operator ${\bf g}$
transforming the antipode map may be, {\it {\`a} la} (\ref{eq:Srel}), 
explicitly evaluated in a closed form:
\beq
{\bf g} = \exp\left(- \frac{1}{2}{\bf TH} (1 - {\bf T}^{-2})\right).
\label{eq:gvalue}
\eeq  
Combining (\ref{eq:gvalue}) with the property (\ref{eq:antwi}) we now 
immediately obtain a disentanglement relation, which, if expressed 
in terms of classical generators, reads as follows:  
\beq
\mu \left[\exp\left(\frac{1}{2}\,h \otimes 
\ln (1 - 2{\bf h} b_{+})\right)\right] = \exp(-{\bf h} h b_{+}).
\label{eq:disentan}
\eeq
In the above relation ${\bf h}$ may be treated as an arbitrary
parameter. To our knowledge the above disentanglement formula 
involving the classical $sl(2)$ generators was not observed before. 

\par

Another map, where the Cartan element ${\bf H}$ of the deformed
$\Uh$ algebra remains diagonal, is given by 'mapping functions' 
\beq
\varphi_1(b_{+}) = \frac{1}{\sqrt{1 - \frac{{\bf h}^2 b_{+}^2}{4}}},
\qquad
\psi_1({\bf T}) = \hbox{sech}(\frac{\bf hX}{2}).
\label{eq:2ndcho}
\eeq 
The $\Uh$ algebra may now be mapped on the the classical $\U$ algebra
as 
\bea
{\bf E} &=& \frac{1}{\sqrt{1 - \frac{{\bf h}^2 b_{+}^2}{4}}}\,e,\qquad
{\bf H} = h,\qquad {\bf T}^{\pm 1} 
= \frac{1 \pm \frac {{\bf h} b_{+}}{2}}{1 \mp \frac {{\bf h} b_{+}}{2}},  
\nn\\
{\bf F} &=& \sqrt{1 - \frac{{\bf h}^2 b_{+}^2}{4}}\,f
- \frac{{\bf h}^2 b_{+}}
{4 \left(1 - \frac{{\bf h}^2 b_{+}^2}{4}\right)^{\frac{3}{2}}}\,e
- \frac{{\bf h}^2 b_{+}}
{4\,\sqrt{1 - \frac{{\bf h}^2 b_{+}^2}{4}}}\,eh.
\label{eq:hdiadir}
\eea
The inverse map now reads
\bea
e &=& \hbox{sech}(\frac{\bf hX}{2}) {\bf E},\qquad\qquad h = {\bf H}, 
\nn\\
f &=& \hbox{cosh}(\frac{\bf hX}{2})\,{\bf F} 
+ \frac{\bf h}{4} \hbox{sinh}({\bf hX})\,
\hbox{cosh}(\frac{\bf hX}{2})\,{\bf E} 
+ \frac{\bf h}{2} \hbox{sinh}(\frac{\bf hX}{2})\,{\bf EH}. 
\label{eq:hdiainv}
\eea
The twist operator ${\cal G}$ for the map (\ref{eq:hdiadir}), unlike the 
previous example of closed-form expression for the 'minimal twist map'
given in (\ref{eq:mtmtwist}), may be determined only in a series:
\beq
{\cal G} = {\bf I} \otimes {\bf I} + \frac{\bf h}{2} {\bf r}
+ \frac{{\bf h}^2}{8} ({\bf r}^2 + {\bf H} \otimes {\bf X}^2 
+{\bf X}^2 \otimes {\bf H}) + O({\bf h}^3),
\label{eq:hdiatwi}
\eeq
where ${\bf r} = {\bf H} \otimes {\bf X} - {\bf X} \otimes {\bf H}$.
The corresponding transforming operator for the antipode map may be 
easily determined from (\ref{eq:Srel}):
\beq
{\bf g} = 1 - {\bf h X} + \frac{1}{2} {\bf h}^2 {\bf X}^2 
+ O({\bf h}^3).
\label{eq:hdiant}
\eeq 

\setcounter{equation}{0}
\section{Conclusion}

In the present work we have studied two aspects of the deformations
of the $osp(2|1)$ superalgebra. Generalizing the approach in
Ref. \cite{K98} we, in Section 2, have found a generic 
representation-independent way of extracting various structures like 
arbitrary finite dimensional $R_{\sf h}$ matrices, $L$ operators and 
${\cal T}$ matrices for the $\cUh$ algebra - obtained {\it via} 
quantization of the $r_2$ matrix - from the corresponding quantities 
of the standard $q$-deformed $\Uq$ algebra. This approach may be 
used, for instance, to construct the higher dimensional ${\cal T}$ 
matrix elements and the corresponding non-commutative spaces
invariant under the coaction of these ${\cal T}$ matrix elements.  
This will be elaborated elsewhere.
\par

In Section 3 of the present work we have recast the algebra  
$\Uh$ algebra - generated by the quantization of the classical $r_1$
matrix - in terms of the nonlinear basis elements. There are several 
benefits of doing this. Ohn's Jordanian deformation \cite{O92} of 
the $sl(2)$ algebra explicitly appear as a Hopf subalgebra
of our $\Uh$ algebra. Moreover, our coproduct rules are considerably
simple. Our approach is expected to be useful in constructing, as
advocated in \cite{BH99} physical many-body models of deformed 
integrable systems obeying the coalgebra symmetry. Explicit
construction of such models based on the coproduct structures of the
$osp(2|1)$ algebras will be presented elsewhere. Of particular interest 
are the possible physical effects arising out of the distinct coproduct
structures of the two quantized algebras $\cUh$ and $\Uh$ discussed in 
Sections 2 and 3 respectively.     

\vskip 0.5cm 
\noindent {\bf Acknowledgments:} 
 
We thank N. Aizawa, A. Chakrabarti and J. Segar for fruitful discussions. 
One of us (B.A.) wants to thank Abir and Ahmed-Nour for their precious 
help. Part of the work was when R.C. visited Institute of Mathematical 
Sciences, Chennai, 600 113,India. He is partially supported by the grant
DAE/2001/37/12/BRNS, Government of India.

\end{document}